# Convexity and sandwich theorems

## Flavia-Corina MITROI-SYMEONIDIS

Lumina – The University of South-East Europe, Faculty of Engineering Sciences, Șos. Colentina 64B, RO-021187, Bucharest, Romania

**Email address:**
fcmitroi@yahoo.com

**Abstract:** We review sandwich theorems from the theory of convex functions.

## 1. Convexity and set-valued functions revisited

Let I be an open interval. The convexity of a function $f$: I $\to \mathbb{R}$ means that it holds
$$tf(x)+(1-t)f(y) \geqslant f(tx+(1-t)y),$$
for all $x, y \in$ I, $t \in [0,1]$.

Recently, convexity has been the subject of intensive research. In particular, many improvements, generalizations and applications of it can be found in the literature.

We denote by $n(\mathbb{R})$ the family of all non-empty subsets of $\mathbb{R}$ and by $cl(\mathbb{R})$ the family of all non-empty and closed subsets of $\mathbb{R}$. A set-valued function F:I$\to n(\mathbb{R})$ is said to be *convex* if it satisfies
$$tF(x)+(1-t)F(y) \subset F(tx+(1-t)y),$$
for all $t \in [0,1]$ and $x, y$ from its domain.

The notions of *concave* and *affine set-valued* function are also considered, when
$$tF(x)+(1-t)F(y) \supset F(tx+(1-t)y),$$
respectively when the two sets coincide for all $t \in [0,1]$ and $x, y$ from the domain of definition. See also [8].

It has been proved in [7] (see also [10]) the following "sandwich" result:

**Theorem 1** Let I be an interval and $f, g$: I$\to \mathbb{R}$. Then the following conditions are equivalent:

i) there exists an affine function $h$:I$\to \mathbb{R}$ such that
$$f(x) \leqslant h(x) \leqslant g(x)$$
on I;

ii) there exists a convex function $h_1$:I$\to \mathbb{R}$ and a concave one $h_2$:I$\to \mathbb{R}$ such that
$$f(x) \leqslant h_1(x) \leqslant g(x) \text{ and } f(x) \leqslant h_2(x) \leqslant g(x)$$
on I;

iii) for all $x, y \in$ I and $t \in [0,1]$,
$$f(tx+(1-t)y) \leqslant tg(x)+(1-t)g(y)$$
and
$$g(tx+(1-t)y) \geqslant tf(x)+(1-t)f(y).$$

For more details about the convex functions see for instance the monograph of C. P. Niculescu and L.-E. Persson [6].

A counterpart of this theorem in the framework of set-valued functions has been recently proved by the author [3]:

**Theorem 2** Let I be an open interval. Let F, G :I$\to cl(\mathbb{R})$ be two set-valued functions. Then the following statements are mutually equivalent:

i) there exists an affine set-valued function H:I$\to cl(\mathbb{R})$ such that
$$F(x) \supset H(x) \supset G(x)$$
on I;

ii) there exists a convex set-valued function $H_1$:I$\to cl(\mathbb{R})$ and a concave set-valued function $H_2$:I$\to cl(\mathbb{R})$ such that
$$F(x) \supset H_1(x) \supset G(x) \text{ and } F(x) \supset H_2(x) \supset G(x)$$
on I;

iii) the functions F and G satisfy
$$F(tx+(1-t)y) \supset tG(x)+(1-t)G(y)$$
and
$$G(tx+(1-t)y) \subset tF(x)+(1-t)F(y).$$

It is known [2] that if F :I$\to cl(\mathbb{R})$ is a convex set-valued function then it has one of the following forms:
  a) $F(x)=[f_1(x), f_2(x)]$
  b) $F(x)=[f_1(x), \infty)$
  c) $F(x)=(-\infty, f_2(x)]$
  d) $F(x)=\mathbb{R}$.

Here $f_1$:I$\to \mathbb{R}$ is a convex function and $f_2$:I$\to \mathbb{R}$ is a concave function.

## 2. Alternative proof of a convexity result

We now provide a simpler proof of Lemma 2 in [5]. For the reader's convenience, we insert here the statement of it:

**Proposition 3** Let $\phi$ and $\psi$ be two functions on an interval I

such that $\psi-\phi$ is increasing (resp. decreasing) on I and $\psi$ is convex (resp. concave) on I. Then
$$(1-t)\phi(x)+t\psi(y) \geq ((1-t)\phi+t\psi)((1-t)x+ty) \text{ (resp } \leq),$$
for all $t \in (0,1)$ and all $x,y \in I$, $x \leq y$.

*Proof*

Let $\psi-\phi$ be increasing and $\psi$ convex. Mutatis mutandis, the other case can be proved similarly.

Due to the monotonicity assumption one has
$$\psi((1-t)x+ty)-\psi(x) \geq \phi((1-t)x+ty)-\phi(x)$$
Using (1) and the convexity of $\psi$ on I, we obtain
$$(1-t)\phi(x)+t\psi(y)-((1-t)\phi+t\psi)((1-t)x+ty)$$
$$= t(\psi(y)-\psi((1-t)x+ty))-(1-t)(\phi((1-t)x+ty)-\phi(x))$$
$$\geq t(\psi(y)-\psi((1-t)x+ty))-(1-t)(\psi((1-t)x+ty)-\psi(x))$$
$$= t\psi(y)+(1-t)\psi(x)-\psi((1-t)x+ty) \geq 0$$
for all $t \in (0,1)$. The proof is completed.

The first of these inequalities holds with equality sign if and only if $\psi-\phi$ is constant and the last one if and only if $\psi$ is affine.

Notice that the particular case $\psi=\phi$ satisfies the hypothesis in Proposition 3, but then the conclusion just degenerates to the definition of a convex (resp. concave) function.

**Open problem** Is there any counterpart of this result in the framework of convex set-valued functions?

## 3. Popoviciu's inequality revisited

Fifty years ago Tiberiu Popoviciu published the following characterization of convex functions [9]:

"A real-valued continuous function $f$ defined on an interval I is convex if and only if it verifies the inequality
$$\frac{f(x)+f(y)+f(z)}{3}+f\left(\frac{x+y+z}{3}\right)$$
$$\geq \frac{2}{3}\left(f\left(\frac{x+y}{2}\right)+f\left(\frac{y+z}{2}\right)+f\left(\frac{x+z}{2}\right)\right)$$
whenever $x,y,z \in I$."

For set-valued functions we see that:

**Proposition 4** A convex set-valued continuous function $F: I \to cl(\mathbb{R})$ verifies the inclusion
$$\frac{F(x)+F(y)+F(z)}{3}+F\left(\frac{x+y+z}{3}\right)$$
$$\subset \frac{2}{3}\left(F\left(\frac{x+y}{2}\right)+F\left(\frac{y+z}{2}\right)+F\left(\frac{x+z}{2}\right)\right)$$
whenever $x,y,z \in I$.

The converse also holds true, via an analogous reasoning as for Popoviciu's real-valued case, since if $F: I \to cl(\mathbb{R})$ is a continuous set-valued function of the form $F(x)=[f_1(x),f_2(x)]$ for all $x \in I$, then the functions $f_1$ and $f_2$ are continuous.

We considered the notion of *continuous set-valued function* according to [2]: A set-valued function $F: I \to n(\mathbb{R})$ is said to be *continuous at a point* $x_0 \in I$ if for every neighborhood V of zero there exists a neighborhood U of zero such that $F(x) \subset F(x_0)+V$ and $F(x_0) \subset F(x)+V$ for all $x \in (x_0+U) \cap I$.

In [1] we find the following lemma:

**Lemma 5** Let $f:[a,b] \to \mathbb{R}$ be a convex function. If $x_1,\ldots,x_n \in [a,b]$ and a convex combination $\sum_{i=1}^n \mu_i x_i$ of these points equals a convex combination $\lambda_1 a + \lambda_2 b$ of the endpoints, then
$$\sum_{i=1}^n \mu_i f(x_i) \leq \lambda_1 f(a) + \lambda_2 f(b).$$

Hence we notice that if we consider the particular case $x_1=\frac{x+y}{2}$, $x_2=\frac{y+z}{2}$, $x_3=\frac{x+z}{2}$ with equal weights $\mu_i=\frac{1}{3}$ for $i=1,2,3$, then we also find another upper bound of the right hand side term of Popoviciu's inequality:
$$\frac{2}{3}\left(f\left(\frac{x+y}{2}\right)+f\left(\frac{y+z}{2}\right)+f\left(\frac{x+z}{2}\right)\right) \leq f(a)+f(b)$$
whenever one has $a$ and $b$ such that $x,y,z \in [a,b]$ and
$$\frac{x+y+z}{3} = \frac{a+b}{2}.$$

Moreover, under these conditions we get as particular cases of Lemma 5 the inequalities
$$2\frac{f(x)+f(y)+f(z)}{3} \leq f(a)+f(b)$$
and
$$2f\left(\frac{x+y+z}{3}\right) \leq f(a)+f(b).$$

By summing the above inequalities, we obtain the following statement:

**Proposition 6** A real-valued continuous convex function $f$ defined on an interval $[a,b]$ verifies the double inequality
$$f(a)+f(b) \geq \frac{f(x)+f(y)+f(z)}{3}+f\left(\frac{x+y+z}{3}\right)$$
$$\geq \frac{2}{3}\left(f\left(\frac{x+y}{2}\right)+f\left(\frac{y+z}{2}\right)+f\left(\frac{x+z}{2}\right)\right)$$
for all $x,y,z \in [a,b]$ such that $\frac{x+y+z}{3} = \frac{a+b}{2}$.

We will now establish a corresponding version of Lemma 5 for set-valued functions:

**Proposition 7** Let $F:[a,b] \to cl(\mathbb{R})$ be a set-valued convex function. If $x_1,\ldots,x_n \in [a,b]$ and a convex combination $\sum_{i=1}^n \mu_i x_i$ of these points equals a convex combination $\lambda_1 a + \lambda_2 b$ of the endpoints, then
$$\sum_{i=1}^n \mu_i F(x_i) \supset \lambda_1 F(a) + \lambda_2 F(b).$$

*Proof*

Straightforward, by considering the above four cases. We only consider the case $F(x)=[f_1(x),f_2(x)]$, $x \in [a,b]$. The remaining cases are dealt similarly.

One has
$$\sum_{i=1}^n \mu_i F(x_i) = \left[\sum_{i=1}^n \mu_i f_1(x_i), \sum_{i=1}^n \mu_i f_2(x_i)\right]$$
and
$$\lambda_1 F(a) + \lambda_2 F(b) =$$
$$= [\lambda_1 f_1(a) + \lambda_2 f_1(b), \lambda_1 f_2(a) + \lambda_2 f_2(b)].$$

We apply Lemma 5 to the convex function $f_1:[a,b] \to \mathbb{R}$ and to the concave function $f_2:[a,b] \to \mathbb{R}$.

Hence
$$\sum_{i=1}^n \mu_i f_1(x_i) \leq \lambda_1 f_1(a) + \lambda_2 f_1(b)$$
and
$$\sum_{i=1}^n \mu_i f_2(x_i) \geq \lambda_1 f_2(a) + \lambda_2 f_2(b)$$
This completes the proof.

For additional recent results connected to the convex set-valued functions the reader is referred to [4].